\documentclass[12pt]{article}

\usepackage[margin=1.25in]{geometry}

\usepackage{amssymb,amsmath}
\usepackage{amsthm}
\usepackage{hyperref}
\usepackage{mathrsfs}
\usepackage{textcomp}
\hypersetup{
    colorlinks,%
    citecolor=red,%
    filecolor=black,%
    linkcolor=blue,%
    urlcolor=black
}

\usepackage{verbatim}
\usepackage{sectsty}
\usepackage{esint}

\usepackage{tikz}

\makeatletter

\newdimen\bibspace
\renewenvironment{thebibliography}[1]{%
 \section*{\refname 
       \@mkboth{\MakeUppercase\refname}{\MakeUppercase\refname}}%
     \list{\@biblabel{\@arabic\c@enumiv}}%
          {\settowidth\labelwidth{\@biblabel{#1}}%
           \leftmargin\labelwidth
           \advance\leftmargin\labelsep
           \itemsep\bibspace
           \parsep\z@skip     %
           \@openbib@code
           \usecounter{enumiv}%
           \let\p@enumiv\@empty
           \renewcommand\theenumiv{\@arabic\c@enumiv}}%
     \sloppy\clubpenalty4000\widowpenalty4000%
     \sfcode`\.\@m}
    {\def\@noitemerr
      {\@latex@warning{Empty `thebibliography' environment}}%
     \endlist}

\makeatother

\makeatletter

\newtheorem{thm}{Theorem}[section]

\newtheorem{prop}[thm]{Proposition}

\newtheorem{rem}[thm]{Remark}

\def\XXint#1#2#3{{\setbox0=\hbox{$#1{#2#3}{\int}$}
  \vcenter{\hbox{$#2#3$}}\kern-.5\wd0}}

\newcommand{\bpf}{\begin{proof}}      \newcommand{\epf}{\end{proof}}
\newcommand{\bpr}{\begin{prop}}      \newcommand{\epr}{\end{prop}}
\newcommand{\brem}{\begin{rem}}      \newcommand{\erem}{\end{rem}}
\newcommand{\bthm}{\begin{thm}}      \newcommand{\ethm}{\end{thm}}

\newcommand{\om}{\Omega}                \newcommand{\pa}{\partial}
           
\newcommand{\be}{\begin{equation}}      \newcommand{\ee}{\end{equation}}

\newcommand{\R}{\mathbb{R}}

\newcommand{\divg}{\mbox{div}}

\begin{document}

\title{\textbf{Green's functions for the heat and Laplace equations with  dynamical boundary conditions in a ball}
\bigskip}

\author{\large  Xuzhou Yang \\
}
\date{}

\maketitle

\begin{abstract}

The Green's functions for the heat and Laplace equations with dynamical boundary conditions in a ball are studied. First, the Green's functions of the Laplace equation with a dynamical boundary condition are given, and the properties of related heat kernels are discussed. Then using these ingredients, two complementary approximations to the heat equation with a dynamical boundary condition in a ball are constructed, including an approximation of Green's function and an approximation of solution. Moreover, the Green's function of the heat equation with a dynamical boundary condition is implicitly presented by a series of eigenfunctions.

\end{abstract}


\smallskip


\bigskip

\tableofcontents

\section{Introduction}\label{sec:intro}

In this paper, we are concerned with some kernels related to the Green's functions to the Laplace equation with a dynamical boundary condition
\be\label{eq:Laplace--dynamical}
  \left\{
\begin{aligned}
 -\Delta u & = 0 \quad                  \mbox{in }   \om \times (0,\infty),\\
 \pa_t u+\pa_\nu u  &    =0    \quad  \mbox{on } \pa \om \times (0,\infty) ,\\
 u(0) &    = \phi_b(x)   \quad  \mbox{on } \pa \om ,
\end{aligned}
\right.
\ee
and the heat equation with a dynamical boundary condition
\be\label{eq:heat--dynamical}
  \left\{
\begin{aligned}
\pa_t u -\Delta u & = 0 \quad                  \mbox{in }   \om \times (0,\infty),\\
 \pa_t u+\pa_\nu u  &    =0    \quad  \mbox{on } \pa \om \times (0,\infty) ,\\
 u(0) &    = \phi_b(x)   \quad  \mbox{on } \pa \om, \\
 u(0) &    = \phi_i(x)   \quad  \mbox{on }  \om,
\end{aligned}
\right.
\ee
for the case that $\om$ is the unit ball of $\R^n$. In the above equations,  $\Delta $ is the Laplace operator with respect to the spatial variable $x\in \om$, $\nu$ is the unit outer normal to $\pa \om$, $\pa_{\nu}$ is the outer normal derivative.

The dynamical boundary condition
\[
\pa_t u + \pa_\nu u = 0 \quad  \mbox{on } \pa \Omega \times (0,\infty)
\]
describes diffusion through the boundary, see \cite{1975Crank}. For example, such diffusion includes the thermal contact with a perfect conductor or diffusion of solute from a well-stirred fluid or vapour.  In particular, it also describes a heat conduction process in $\om$, with a heat source on $\pa \om$, which may depend on the heat flux around the boundary and on the heat flux across it, see \cite{gal2014}.

The Laplace equation with dynamical boundary conditions was already a very famous and classical problem, and many results were established. In \cite{1997fila}, the global solutions of \eqref{eq:Laplace--dynamical} were studied by Fila and Quittner. In \cite{1997amann}, Amann-Fila used such existence conclusions of the linear problem \eqref{eq:Laplace--dynamical} to analyse the Fujita type theorem to related nonlinear occasion of
\eqref{eq:Laplace--dynamical}. They used the theory of $C_0$ semi-group
to generate a Green's function to problem \eqref{eq:Laplace--dynamical} for $\om = \R^n_+$. Then, they applied such representation to the nonlinear problem
\be\label{eq:Laplace--dynamical-fujita}
  \left\{
\begin{aligned}
 -\Delta u & = 0 \quad                  \mbox{in }   \R^n_+ \times (0,\infty),\\
 \pa_t u+\pa_\nu u  &    =u^q    \quad  \mbox{on } \pa \R^n_+ \times (0,\infty) ,\\
 u(0) &    = \phi_b(x)   \quad  \mbox{on } \pa \R^n_+ ,
\end{aligned}
\right.
\ee
for $q \in (1,\infty)$, obtaining a critical exponent for global existence of positive solutions of the above Laplace  equation \eqref{eq:Laplace--dynamical-fujita} in the upper half-space with a nonlinear dynamical boundary condition, which is a Fujita type theorem.

2009, V\'azquez-Vitillaro \cite{vaz2009} studied the well-posedness and regularity of the Laplace equation with dynamical boundary conditions
\be\label{eq:Laplace--dynamicals}
  \left\{
\begin{aligned}
- \Delta u & = 0 \quad                  \mbox{in }   \om \times (0,\infty),\\
 \pa_t u  - \delta \Delta_{\pa \om} u + d \pa_\nu u    &    =0    \quad  \mbox{on } \pa \om \times (0,\infty) ,\\
 u(0) &    = \phi_b(x)   \quad  \mbox{on } \pa \om ,
\end{aligned}
\right.
\ee
in which they concluded that \eqref{eq:Laplace--dynamicals} is ill-posed when $ d <0$ and well-posed when $d \ge 0$ and a similar situation occurs when the Laplacian operator $\Delta$ is replaced by the heat operator $\pa_t - \Delta$ .

2014, Gal-Meyries \cite{gal2014} studied an elliptic equation with dynamical boundary conditions
\be\label{eq:elliptic--dynamicals}
  \left\{
\begin{aligned}
\lambda u  - d \Delta u & = f(u) \quad                  \mbox{in }   \om \times (0,\infty),\\
 \pa_t u - \delta \Delta_{\pa \om} u +d \pa_\nu u  &    =g(u)    \quad  \mbox{on } \pa \om \times (0,\infty) ,\\
 u(0) &    = \phi_b(x)   \quad  \mbox{on } \pa \om ,
\end{aligned}
\right.
\ee
where $d>0 $, $\delta \ge 0$ , $\om$ is a bounded domain with smooth boundary $\pa \om$, $ \Delta_{\pa \om}$ is the Laplace-Beltrami operator
and $\pa_\nu$ is the outer normal derivative on $\pa \om$, $f,g \in C^{\infty}(\R)$ and with assumption that $f$ is
globally Lipschitz continuous and that $\lambda$ is sufficiently large, in dependence on $f$. In \cite{gal2014}, the blow-up, global existence, global attractors and convergence to single equilibria were considered.

Two types of dynamical boundary conditions are buried in above problem \eqref{eq:elliptic--dynamicals}. For $\delta>0$, the dynamical boundary conditions are of reactive-diffusive type, and for $\delta=0$, the dynamical
boundary conditions are purely reactive. The motivation of dynamical boundary conditions comes from physics. The function $u $ represents
the steady state temperature in a body $\om$ such that the rate at which $u$ evolves through the boundary $\pa \om$ is proportional to the flux on the boundary, up to some correction $\delta \Delta_{\pa \om} u (\delta \ge 0)$ which from a modelling viewpoint, accounts for small diffusive effects along $\pa \om$.

In 2024, the nonlinear boundary
 diffusion problem related to \eqref{eq:Laplace--dynamical} was studied by Jin-Xiong-Yang \cite{2024JXY}. Such problem is written as
\begin{equation}\label{eq:BDE}
 \left\{
\begin{aligned}
&  -\Delta u = 0 \quad                           \mbox{in }   \Omega \times (0,\infty),\\
&\partial_{t}(\beta(u))+\partial_\nu u  =  0 \quad    \mbox{on } \pa \Omega  \times  (0,\infty) , \\
&u(x,0)=u_0(x)    \quad      \mbox{on }  \pa \Omega,
\end{aligned}
\right.
\end{equation}
which is closely related to the boundary Yamabe problem, see \cite{brendle2002}. In \cite{2024JXY}, we studied the existence, lower and upper bounds of the positive solution of the nonlinear boundary diffusion problem \eqref{eq:BDE} for $\beta(u)=u^p$ and classified the long time behavior in terms of the exponent $p$.  Especially, we proved the short time existence of smooth positive solution and obtained the long time existence by extending such solution whenever they were positive.  Then we also obtained the stability of the separable solutions of  \eqref{eq:BDE} for both $0<p<1$ and $1<p<\frac{n}{n-2}$ with sharp convergence rates.

Moreover, Capitanelli-D'Ovidio \cite{2024capitanelli} displayed the connection between Brownian motion and equation \eqref{eq:Laplace--dynamical}. They provided a probabilistic formula for the solution of \eqref{eq:Laplace--dynamical}. Therefore, such dynamical boundary condition is linked with probability theory as well as the Dirichlet heat kernel. We may look into such direction in the future research.

As for the parabolic equations with dynamical boundary conditions, there were much research. (see e.g., \cite{bandle2006,below2003,below2000,vaz2011}). Particularly, for heat equation with dynamical boundary conditions, in 2012, Fila-Ishige-Kawakami \cite{2012fila} studied the convergence to Poisson kernel with equation \eqref{eq:heat--dynamical}. In 2020,
Fila-Ishige-Kawakami-Lankeit \cite{2020ex-ball} analysed the connection between solutions of \eqref{eq:Laplace--dynamical} and \eqref{eq:heat--dynamical} for $\Omega$ be the exterior of the unit ball. Then Fila-Ishige-Kawakami \cite{2021fila} studied such connection for $\Omega=\R^n_+$. Above two papers both investigated this connection by looking at an approximation problem
\be\label{eq:heat--dynamical-epsilon}
  \left\{
\begin{aligned}
\epsilon \pa_t u_\epsilon -\Delta u_\epsilon & = 0 \quad                  \mbox{in }   \om \times (0,\infty),\\
 \pa_t u_\epsilon+\pa_\nu u_\epsilon  &    =0    \quad  \mbox{on } \pa \om \times (0,\infty) ,\\
 u_\epsilon(0) &    = \phi_b(x)   \quad  \mbox{on } \pa \om, \\
 u_\epsilon(0) &    = \phi(x)   \quad  \mbox{on }  \om.
\end{aligned}
\right.
\ee
For the two regions that we mentioned above, they proved that the solution $u_\epsilon$ converges to the solution $u$ of \eqref{eq:Laplace--dynamical} as $\epsilon \to 0^+$ and found the upper and lower bounds for the rate of convergence.

More results have been obtained for the case that $\om =\R_+^n$ in terms of the initial boundary problem \eqref{eq:heat--dynamical} recently. Due to the linearity and reflection principle in such case, the explicit formula of the Green's function to \eqref{eq:heat--dynamical} was built recently in 2024 by Ishige-Katayama-Kawakami \cite{2024heat-fundamental} for $\om = \R^n_+$. We briefly exhibit these results in Section \ref{sec:half space} as comparison to the case that $\om = B_1$.

However, for $\om = B_1$, yet there is few result for dynamical boundary problems \eqref{eq:Laplace--dynamical} and \eqref{eq:heat--dynamical}, nor as the Fujita type theorem for the nonlinear case. Therefore, we are curious about the Green's functions for such dynamical boundary problems in a ball.

The nonlinear case of \eqref{eq:heat--dynamical} is called boundary heat control problem. Early in 1972, the physical model of boundary heat control was carried out by Duvaut-Lions \cite{DL1972}. Later in 2010, Athanasopoulos-Caffarelli \cite{2010Continuity} studied the boundary heat control model
\begin{equation}\label{eq:BPME}
 \left\{
\begin{aligned}
& \alpha\partial_{t}u  =\Delta u \quad                           \mbox{in }   \Omega \times (0,T],\\
&-\partial_\nu u  \in  \partial_{t}(\beta(u)) \quad    \mbox{on }   \Gamma \times  (0,T] , \\
&u              =  0  \quad   \mbox{on } \left( \partial \Omega - \Gamma \right) \times (0,T],\\
&u(x,0)=u_0(x)    \quad      \mbox{on }  \Omega.
\end{aligned}
\right.
\end{equation}
They proved the H\"older continuity of the weak solution to \eqref{eq:BPME} for $\beta(u)=u^q,0<q<1$, and the sign ``$\in$" be ``$=$",
which is an approximation to the Stefan problem described in \cite{DL1972}. The mathematical model \eqref{eq:BPME} can be used to describe the temperature in ice-water mixture. Such nonlinear problem in a ball is absent of research as well.

In the present paper, we aim to study the Green's functions and further find some useful kernels for problem \eqref{eq:Laplace--dynamical} and \eqref{eq:heat--dynamical} for $\om = B_1$.  We first give the explicit formula of the Green's function $K_1(x,y,t)$ of the initial value problem of Laplacian equation with a dynamical boundary condition \eqref{eq:Laplace--dynamical} for $\om = B_1$ and display its properties. Then we provide results about Dirichlet heat kernel $\Gamma_1(x,y,t)$ in a ball and built some useful kernels and discuss their properties as preparation, see Proposition \ref{prop:soluton-u-ball-heat kernel-F_1-Gamma_1} and Proposition \ref{prop:soluton-u-ball-heat kernel-H_1}.

Our goal is to employ the above kernels to build approximations for the Green's function $G_1(x,y,t)$ of the initial value problem of heat equation with a dynamical boundary condition \eqref{eq:heat--dynamical} in a ball, for which the implicit representation and related property are discussed in Section \ref{sec:green}. The exploration of explicit formula of $G_1(x,y,t)$ is of much difficulty because we are lack of the reflection principle compared to the case in which $\Omega = \R^n_+$. So, instead we utilize the kernels above to find two approximations for the Green's function and the solution of \eqref{eq:heat--dynamical}, respectively.

Inspired by Ishige-Katayama-Kawakami \cite{2024heat-fundamental}, in which they provided an explicit formula of Green's function of initial value problem \eqref{eq:heat--dynamical} for $\om=\R^n_+$, we build an approximation $\mathcal{G}_1(x,y,t)$ for the Green's function  $G_1(x,y,t)$ of \eqref{eq:heat--dynamical} with $\Omega = B_1 - \{0\}$, such that
\[
|G_1(x,y,t) -\mathcal{G}_1(x,y,t) | \to 0, \quad  t \to 0^+,
\]
for $x,y \in B_1 - \{0\}$, see Theorem \ref{thm:G-1 approximation}.

In \cite{2020ex-ball}, the solution of \eqref{eq:heat--dynamical} is represented  via combinations of solutions to the Laplace equation with an inhomogeneous dynamical boundary condition and an inhomogeneous heat equation with the homogeneous Dirichlet boundary condition for the case that $\Omega $ is the exterior of the unit ball. Follow this lead, using the Green's function $K_1(x,y,t)$ of initial value problem of Laplace equation \eqref{eq:Laplace--dynamical} and the Dirichlet heat kernel $\Gamma_1(x,y,t)$ as ingredients, we build an approximation  to solution of \eqref{eq:heat--dynamical} for $\Omega = B_1 $, as $ t\to 0^+, \  |x| \to 0^+$.  This can be seen as a complement of the approximation solution generated by $\mathcal{G}_1(x,y,t)$, since it provides information at the origin, see Theorem \ref{thm:u approximation}.

This paper is organized as follows. In section \ref{sec:half space}, we gather some recent results for $\om=\R^n_+$ as comparison. In Section \ref{sec:unit ball}, we establish related results for the Dirichlet heat kernel and some auxiliary kernels for $\om=B_1$ as preparation. Section \ref{sec:green} displays the implicit formula and some properties of the Green's function of \eqref{eq:heat--dynamical}. In Section \ref{sec:appro}, we explicitly construct two types of approximations of \eqref{eq:heat--dynamical} by using $K_1$, $\Gamma_1$ and their related kernels.

\section{The case in the upper half space of $\R^n$}\label{sec:half space}

In this section, we recall some notions of Green's functions and exhibit some results of the dynamical boundary problems \eqref{eq:Laplace--dynamical} and \eqref{eq:heat--dynamical} for  $\om=\R^n_+$. For the knowledge of the Green's functions and the Poisson kernels on the half space and unit ball, we refer to the Chapter 2 in the famous monograph of L.C. Evans \cite{evans2010}.
\subsection{The Laplace equation}
Let $G_+(x,y)$ be the Green's function to the Dirichlet Laplacian on the upper half space, taking the form
\[
G_+(x,y)=\Phi_n(y-x)-\Phi_n(y-x^*), \quad x,y\in \R^n_+, x \neq y,
\]
where
\[
\Phi_n(x) =   \left\{
\begin{aligned}
  & -\frac{1}{2\pi} \ln |x| , n=2,    \\
 &    c_n \frac{1}{|x|^{n-2}} , n \ge 3 ,
\end{aligned}
\right.
\]
is the fundamental solution of Laplace equation in $\R^n$ and
\[ x^*= (x_1,...,x_{n-1},-x_n) , \quad  x\in \R^n_+,\]
is the dual point to $x=(x_1,...,x_{n-1},x_n)$ with respect to the hyper plane $\pa \R^n_+$. Then, one can symbolically write
\be\label{eq:Green-Laplace-halfspace-pde}
  \left\{
\begin{aligned}
 -\Delta_y G_+& = \delta_x(y) \quad                           \mbox{in }   \R^n_+,\\
 G_+ &    =0    \quad  \mbox{on } \pa \R^n_+ ,
\end{aligned}
\right.
\ee
where $\delta_x$ denotes the Dirac measure giving unit mass to the point $x$ and $G_+$ considered as a function of $y$.
Then the Poisson kernel on the half space is presented as
\be \label{eq:poisson kernel-halfspace}
P_+(x,y)= -\pa_{\nu_y}G_+(x,y)=c_n \frac{x_n}{|x-y|^n},  \quad x \in \R^n_+, y \in \pa \R^n_+,
\ee
where $\nu_y$ is the outward normal direction at  $y \in \pa \R^n_+$, and $c_n$ is the constant such that
\be\label{eq:poisson-Laplace-plane-integral-P_1}
 \int_{\pa \R^n_+  } P_+(x,y) d y =1 ,    \quad \forall x \in \R^n_+.
\ee
It is symbolically interpreted as
\be\label{eq:poisson-Laplace-halfspace-pde}
  \left\{
\begin{aligned}
 -\Delta_x P_+ & = 0 \quad                  \mbox{in }   \R^n_+,\\
 P_+ &    = \delta_x(y)   \quad  \mbox{on } \pa \R^n_+ ,
\end{aligned}
\right.
\ee
with $P_+$ considered as a function of $x$.

Let

\[
K_+(x,y,t)=P_+(x+te_n,y), \quad x\in \R^n_+, y\in \pa \R^n_+,t>0,
\]
for $e_n$ denoting the unit vector $(0,...,0,1) \in \R^n$. By the property of the Poisson kernel $P_+(x,y)$, one concludes that $K_+(x,y,t)$ is the Green's function to \eqref{eq:Laplace--dynamical} for $\Omega = \R^n_+$, which is symbolically interpreted as
\be\label{eq:poisson-Laplace-halfspcace-pde-dynamical}
  \left\{
\begin{aligned}
 -\Delta_x K_+ & = 0 \quad                  \mbox{in }   \R^n_+ \times (0,\infty),\\
 \pa_t K_++\pa_{\nu_x} K_+  &    =0    \quad  \mbox{on } \pa \R^n_+ \times (0,\infty) ,\\
 K_+(0) &    = \delta_x(y)   \quad  \mbox{on } \pa \R^n_+ \times \{0\}, 
\end{aligned}
\right.
\ee
with $K_+$ considered as a function of $(x,t)$ to satisfy above expression \eqref{eq:poisson-Laplace-halfspcace-pde-dynamical}. So
\[
v(x,t)= \int_{\pa \R^n_+} K_+(x,y,t)\phi_b(y)d y
\]
is a solution to \eqref{eq:Laplace--dynamical} and $K_+(x,y,t)$ would satisfy
\be\label{eq:dynamical-Laplace-half space-integral-K_+}
 \int_{\pa R^n_+} K_+(x,y,t) d y =1 ,    \quad \forall x \in R^n_+, t>0.
\ee

\subsection{The Heat equation}

Let $\Gamma_+(x,y,t)$ be the Dirichlet heat kernel on the upper half space, symbolically written as
\be\label{eq:dirichlet-heat-kernel-ball-pde-Gamma_+}
  \left\{
\begin{aligned}
 \pa_t \Gamma_+ -\Delta_x \Gamma_+ & = \delta_x(y)\delta_0(t) \quad                           \mbox{in }   \R^n_+ \times (0, \infty),\\
 \Gamma_+ &    =0    \quad  \mbox{on } \pa \R^n_+\times (0, \infty) ,\\
 \Gamma_+(0) & = \delta_x(y)   \quad  \mbox{on } \R^n_+\times \{0 \}.
\end{aligned}
\right.
\ee
\[
\Gamma_+(x,y,t) = \Gamma(x,y,t) - \phi_+^y(x,t), \quad  x,y \in \R^n_+, t>0,
 \]
where
$\delta_0(t)$ denoting the Dirac measure on $(0,\infty)$ giving unit mass to point $0$, 
\[
\Gamma(x,y,t)= \frac{1}{(4\pi t)^{\frac{n}{2}}} \exp (-\frac{|x-y|^2}{4t} ) , \quad x,y \in \R^n, t>0,
\]
 is the standard heat kernel in the $n$ dimensional Euclidean space and the correction function $\phi_+^y(x,t)$ satisfies
 \be\label{eq:dirichlet-heat kernel-halfspace-correction-pde}
  \left\{
\begin{aligned}
 \pa_t \phi_+^y(x,t) -\Delta_x \phi_+^y(x,t) & = 0 \quad                           \mbox{in }   \R^n_+ \times (0, \infty),\\
 \phi_+^y(x,t) &    =\Gamma(x,y,t)    \quad  \mbox{on } \pa \R^n_+ \times (0, \infty) ,\\
 \phi_+^y(x,t) & = 0 \quad  \mbox{on } \R^n_+ \times \{0 \}.
\end{aligned}
\right.
\ee
Thanks to the reflection principle, the correction function $\phi_+^y(x,t)$ is solved as
\[
\phi_+^y(x,t) = \Gamma(x,y^*,t), \quad  x,y \in \R^n_+, t>0,
\]
which provides us an explicit formula of $\Gamma_+(x,y,t)$:
\be
\begin{split}
\Gamma_+(x,y,t)&=     \Gamma(x,y,t) -       \Gamma(x,y^*,t)\\
&= (4\pi t)^{-\frac{n}{2}} \left[ \exp\left(-\frac{|x - y|^2}{4t}\right) - \exp\left(-\frac{|x - y^*|^2}{4t}\right) \right]\\
&= \Phi_{n-1}(x' - y', t) \left( \Phi_1(x_n - y_n, t) - \Phi_1(x_n + y_n, t) \right)
\end{split}
\ee
where $x=(x',x_n)$, $y=(y',y_n) \in \R^n_+$ , $t>0$ and $\Phi_n(x,t) =\frac{1}{(4\pi t)^{\frac{n}{2}}} \exp (-\frac{|x|^2}{4t} )   $ denotes the fundamental solution of the heat equation in $\R^n $. For $\Omega = \R_+^n$, $\pa_{\nu_x} = -\pa_{x_n}$, equation \eqref{eq:heat--dynamical} is of constant coefficients. In \cite{2024heat-fundamental}, 2024, Ishige-Katayama-Kawakami constructed the Green's function $G_+(x,y,t)$ for the initial value problem \eqref{eq:heat--dynamical}, that
\[
   u(x,t) = \int_{\R_+^n} G_+(x,y,t) \phi_i(y) dy + \int_{\partial \R_+^n} G_+(x,y,t) \phi_b(y) dy
\]
is a classical solution to \eqref{eq:heat--dynamical}, and
\be \label{eq:-heat-ball-integral-G_+}
\int_{\pa \R_+^n}  G_+(x,y,t) d y + \int_{ \R_+^n} G_+(x,y,t) dy = 1, \quad \forall
x \in \R_+^n, t>0,
\ee
\[
 G_+(x,y,t) = G_+(y,x,t) , \quad x,y \in \R_+^n, t>0.
\]
They provided the explicit formula of $G_+(x,y,t)$ by
\be
\begin{split}
G_+(x,y,t)&= \Gamma_+(x,y,t) +\int_0^t \pa_{x_n}\Gamma_+(x+se_n,y,t-s)ds\\
&= \Gamma_+(x,y,t) +\int_0^t -2  \pa_{x_n}\Gamma(x+se_n,y*,t-s)ds.
\end{split}
\ee
Follow this pioneer work, we wish to study the Green's function of \eqref{eq:heat--dynamical} for $\Omega = B_1$.

\section{The case in a ball of $\R^n$} \label{sec:unit ball}

\subsection{The Laplace equation}

Let $G_1(x,y)$ be the Green's function to the Dirichlet Laplacian on the unit ball, taking the form
\[
G_1(x,y)=\Phi_n(y-x)-\Phi_n(|x|(y-\tilde{x})), \quad x,y\in B_1, x \neq y,
\]
where
\[
\Phi_n(x) =   \left\{
\begin{aligned}
  & -\frac{1}{2\pi} \ln |x| , n=2,    \\
 &    c_n \frac{1}{|x|^{n-2}} , n \ge 3 ,
\end{aligned}
\right.
\]
is the fundamental solution of Laplace equation in $n$ dimensional Euclidean space $\R^n$ and
\[ \tilde{x}= \frac{x}{|x|^2} , \quad  x\in \R^n-\{0\},\]
is the dual point to $x$ with respect to the unit sphere $\pa B_1$.
It is symbolically interpreted as
\be\label{eq:Green-Laplace-ball-pde}
  \left\{
\begin{aligned}
 -\Delta_y G_1 & = \delta_x(y) \quad                           \mbox{in }   B_1,\\
 G_1 &    =0    \quad  \mbox{on } \pa B_1 ,
\end{aligned}
\right.
\ee
where $\delta_x(y)$ denotes the Dirac measure giving unit mass to the point $x$ and $G_1(x,y)$ considered as a function of$ y$.
Then the Poisson kernel in the unit ball is presented as
\be \label{eq:poisson kernel-ball}
P_1(x,y)= -\pa_{\nu_y}G_1(x,y)=c_n \frac{1-|x|^2}{|x-y|^n},  \quad x \in B_1, y \in \pa B_1,
\ee
where $\nu_y$ is the outward normal direction at  $y \in \pa B_1$. Here, $c_n$ is the constant such that
\be\label{eq:poisson-Laplace-sphere-integral-P_1}
 \int_{\pa B_1} P_1(x,y) d \sigma_y =1 ,    \quad \forall x \in B_1,
\ee
where $d\sigma_y$ represents the element surface area at $y \in \pa B_1$.  It is symbolically interpreted as
\be\label{eq:poisson-Laplace-ball-pde}
  \left\{
\begin{aligned}
 -\Delta_x P_1 & = 0 \quad                  \mbox{in }   B_1,\\
 P_1 &    = \delta_x(y)   \quad  \mbox{on } \pa B_1 ,
\end{aligned}
\right.
\ee
with $P_1(x,y)$ considered as a function of $x$ to satisfy \eqref{eq:poisson-Laplace-ball-pde}.
\subsubsection{A decomposition using the Green's function to the Dirichlet problem}

Under this scheme we can decompose a function $f(x) \in C^2(B_1) \cap C(\bar{B_1})$ into
\be\label{eq:Laplace-ball-decomposition}
 f(x)=\int_{B_1}G_1(x,y)(-\Delta f (y)) dy + \int_{\pa B_1}P_1(x,y)f(y)d\sigma_y.
\ee

Let
\[
u(x,t)=f(xe^{-t}),
\]
then $u$ satisfies
\be\label{eq:Laplace-ball-dynamical-nonhomogeneous interior}
  \left\{
\begin{aligned}
 -\Delta u & =  -e^{-2t}\Delta f(xe^{-t}) \quad                           \mbox{in }   B_1 \times (0,\infty),\\
 \pa_t u+\pa_{\nu_x u}  &    =0    \quad  \mbox{on } \pa B_1 \times (0,\infty) ,\\
 u(x,0)&   =f(x).
\end{aligned}
\right.
\ee

\subsubsection{Related kernels}
Let
\[
J_1(x,y,t)=G_1(xe^{-t},y),\quad x,y \in B_1, t>0,
\]
and
\[
K_1(x,y,t)=P_1(xe^{-t},y), \quad x\in B_1, y\in \pa B_1,t>0.
\]

By the property of the Poisson kernel $P_1(x,y)$, one concludes that $K_1(x,y,t)$ is the Green's function to \eqref{eq:Laplace--dynamical} for $\Omega=B_1$, which is interpreted as
\be\label{eq:poisson-Laplace-ball-pde-dynamical}
  \left\{
\begin{aligned}
 -\Delta_x K_1 & = 0 \quad                  \mbox{in }   B_1 \times (0,\infty),\\
 \pa_t K_1+\pa_{\nu_x} K_1  &    =0    \quad  \mbox{on } \pa B_1 \times (0,\infty) ,\\
 K_1(0) &    = \delta_x(y )   \quad  \mbox{on } \pa B_1 \times \{0\}.
\end{aligned}
\right.
\ee
So
\[
v(x,t)= \int_{\pa B_1} K_1(x,y,t)\phi_b(y)d\sigma_y
\]
is a solution to \eqref{eq:Laplace--dynamical} and $K_1(x,y,t)$ would satisfy
\be\label{eq:dynamical-Laplace-sphere-integral-K-1}
 \int_{\pa B_1} K_1(x,y,t) d \sigma_y =1 ,    \quad \forall x \in B_1, t>0.
\ee

Secondly, for
\[
w(x,t)= \int_{ B_1} J_1(x,y,t)(-\Delta\phi_i(y))dy
\]
we see that it is a solution to the initial value problem with a dynamical boundary
\be\label{eq:Laplace-ball-dynamical-interior-initial}
  \left\{
\begin{aligned}
 -\Delta w & =  -e^{-2t}\Delta \phi_i(xe^{-t}) \quad                           \mbox{in }   B_1 \times (0,\infty),\\
 \pa_t w+\pa_\nu w  &    =0    \quad  \mbox{on } \pa B_1 \times (0,\infty) ,\\
 w(x,0)&   =\Psi(x) \quad  \mbox{on }  \bar{B_1}  ,
\end{aligned}
\right.
\ee
where the initial data $\Psi(x)$ is a solution to
\be\label{eq:Laplace-ball-initial-data-dirichlet}
  \left\{
\begin{aligned}
 -\Delta \Psi & =  -\Delta \phi_i \quad                           \mbox{in }   B_1 \\
 \Psi  &    =0    \quad  \mbox{on } \pa B_1 .
\end{aligned}
\right.
\ee
$\Psi(x)$ coincides with $\phi_i(x)$, if $\phi_i(x)$ is of class $C^2_c(B_1)$.

\subsection{The Heat equation}
Let $\Gamma_1(x,y,t)$ be the Dirichlet heat kernel on the unit ball, interpreted as
\be\label{eq:dirichlet-heat-kernel-ball-pde-Gamma_1}
  \left\{
\begin{aligned}
 \pa_t \Gamma_1 -\Delta_x \Gamma_1 & = \delta_x(y)\delta_0(t) \quad                           \mbox{in }   B_1 \times (0, \infty),\\
 \Gamma_1 &    =0    \quad  \mbox{on } \pa B_1 \times (0, \infty) ,\\
 \Gamma_1(0) & = \delta_x(y)   \quad  \mbox{on } B_1 \times \{0\}.
\end{aligned}
\right.
\ee
\[
\Gamma_1(x,y,t) = \Gamma(x,y,t) - \phi_1^y(x,t), \quad  x,y \in B_1, t>0,
 \]
where
\[
\Gamma(x,y,t)= \frac{1}{(4\pi t)^{\frac{n}{2}}} \exp (-\frac{|x-y|^2}{4t} ) , \quad x,y \in \R^n, t>0,
\]
 is the standard heat kernel in the $n$ dimensional Euclidean space and the correction function $\phi_1^y(x,t)$ satisfies the boundary value problem
 \be\label{eq:dirichlet-heat kernel-ball-correction-pde}
  \left\{
\begin{aligned}
 \pa_t \phi_1^y(x,t) -\Delta_x \phi_1^y(x,t) & = 0 \quad                           \mbox{in }   B_1 \times (0, \infty),\\
 \phi_1^y(x,t) &    =\Gamma(x,y,t)    \quad  \mbox{on } \pa B_1 \times (0, \infty) ,\\
 \phi_1^y(x,t) & = 0 \quad  \mbox{on } B_1.
\end{aligned}
\right.
\ee

It is noticed that one cannot solve out the expression of such corrector $\phi_1^y(x,t)$, since the  spherical reflection is not so perfect as the planar reflection. So yet, we do not know the exact expression of $\Gamma_1(x,y,t)$. Here, we summarize two types of expressions of $\Gamma_1$. The first one is a probability expression, which is the famous Hunt formula. The second one is an analytic expression, which comes from a series of oscillating components.

\subsubsection{Probability view}
Unlike the Dirichlet heat kernel on the half space, there is no explicit expression for $\Gamma_1(x,y,t)$ because we are lack of the reflection principle. Thanks to the theory of Brownian motion, there are some formulas for Dirichlet heat kernel using probability notations.

We consider n-dimensional Brownian motion $W=(W_t)_{t\ge 0}$ starting from $x \in \R^n$ and we denote by $\mathbf{P}^x$ and $\mathbf{E}^x$ the corresponding probability law and the expected value. $\mathbf{P}^x$ is absolutely continuous with respect to the Lebesgue measure and $\Gamma(x,y,t)$ is the corresponding transition probability density.

For the unit ball $B_1 \subset \R^n$, we define the first exit time from $B_1$ by
\[
\tau_1 = \inf \left\{t > 0 : W_t \notin B_1 \right\}.
\]
Then $\Gamma_1(x,y,t)$ is the transition probability density for $W^1 = (W^1_
t )_{t\ge0}$ Brownian
motion killed upon leaving the unit ball $B_1$. The relation between $\Gamma_1(x,y,t)$ and $\Gamma(x,y,t)$
together with the joint distribution of $(\tau_1, W_{\tau_1} )$ is described by the Hunt formula
\[
\Gamma_1(x,y,t) = \Gamma(x,y,t) - \mathbf{E}^x \left( t>\tau_1, \Gamma(W_{\tau_1} , y,t - \tau_1)\right), \quad x, y \in B_1, t > 0.
\]
Many research of the upper and lower estimates of $\Gamma_1(x,y,t)$ had been done following this probability view. For recent result, see 2020 \cite{2020dirichlet-heatkernel-ball}, in which Maecki-Serafin described the exponential behaviour of $ \Gamma_1(x,y,t)$ when $t$ is small:
\be \label{eq:Gamma-1-estimate}
  \frac{1}{C}h(x,y,t)\Gamma(x,y,t) \le \Gamma_1(x,y,t) \le   C h(x,y,t)\Gamma(x,y,t) , \quad x,y \in B_1, t\in(0,T)
\ee
where
\be  \label{eq:h}
\begin{split}
h(x,y,t) = &\left(     1 \land   \frac{(1-|x|)(1-|y|)}{t}           \right)+\left(        1 \land  \frac{(1-|x|) |x-y|^2}{t}                     \right)\left(        1 \land  \frac{(1-|y|) |x-y|^2}{t}                     \right).
\end{split}
\ee
 Further, they used result \eqref{eq:Gamma-1-estimate} to find
\be\label{eq:E-1-estimate}
\frac{1}{C}l(x,y,t)\Gamma(x,y,t) \le -\pa_{\nu_y}\Gamma_1(x,y,t) \le   Cl(x,y,t)\Gamma(x,y,t) , \quad x\in B_1,y \in \pa B_1, t\in(0,T)
\ee
where
\be \label{eq:l}
l(x,y,t) = \frac{1-|x|}{t} + \frac{|x-y|^2}{t} \left(   1 \land \frac{(1-|x|)|x-y|^2}{t}                 \right).
\ee
In above expressions, $C$ maybe different in different lines, is a positive constant only depending on the dimension $n$ and time $T$ , and
\[
x \land y = \min\{ x,y   \}.
\]
\subsubsection{Analytic view}
On a more general domain $\om$ in $\R^n$, such an explicit formula like the classical formula of $\Gamma(x,y,t)$ is not generally possible. The next simplest cases of a disc or square involve, respectively, Bessel functions and Jacobi theta functions. Nevertheless, the heat kernel still exists and is smooth for $t > 0$ on arbitrary domains and indeed on any Riemannian manifold with boundary, provided the boundary is sufficiently regular. More precisely, in these more general domains, the Dirichlet heat kernel is the solution of the initial boundary value problem

 \be
{\begin{cases}{\frac {\partial \Gamma_\om}{\partial t}}(x,y,t)=\Delta _{x}\Gamma_\om(x,y,t)&{\text{for all }}t>0{\text{ and }}x,y\in \Omega \\[6pt]\lim _{t\to 0}\Gamma_\om(x,y,t)=\delta _{x}(y)&{\text{for all }}x,y\in \Omega \\[6pt]\Gamma_\om(x,y,t)=0& \text{for all } t>0 \  , \  x  {\text{ or }}y\in \partial \Omega \end{cases}}
\ee

Therefore, for $\om = B_1$, we only have an implicit formula of $\Gamma_1(x,y,t)$, by using a series of eigenfunctions.

\[
\Gamma_1(x,y,t) = \sum_{i=1}^{\infty} \exp( -\lambda_i t) \phi_i(x)\phi_i(y), \quad x, y \in B_1, t > 0,
\]
where $(\lambda_i , \phi_i )$ is the eigen pair for the Dirichlet problem
\be\label{eq:dirichlet-laplace-ball-eigen}
  \left\{
\begin{aligned}
 -\Delta \phi_i & = \lambda_i \phi_i \quad                           \mbox{in }   B_1 ,\\
 \phi_i & = 0 \quad  \mbox{on } \pa B_1,
\end{aligned}
\right.
\ee
satisfying
\[
0<\lambda_1 < \lambda_2 \le \lambda_3 \le \cdots, \quad \lambda_i \to \infty
\]
and $\{ \phi_i \}_{i=1}^\infty$ is an orthonormal basis in $L^2(B_1)$ .
Such perspective reveals

\begin{itemize}
      \item   $\Gamma_1(x,y,t) = \Gamma_1(y,x,t) , \quad x, y \in B_1, t > 0$
      \item   $\Gamma_1(x,y,t) \sim (1-|x|)(1-|y|)\exp(-\lambda_1 t) , \quad$ for large $t$.
\end{itemize}

In one word, so far we know that $\Gamma_1(x,y,t)$ behaves like $\Gamma(x,y,t)$ when $x,y$ is away from $\pa B_1$ and $t$ is small, and $\Gamma_1(x,y,t)$ behaves like behaves like $(1-|x|)(1-|y|)\exp(-\lambda_1 t)$ when $t$ is large, see \cite{EBD1987} and \cite{EBD1984}. The boundary behaviour of $\Gamma_1(x,y,t)$ for small $t$, is still a mystery.

\subsubsection{A decomposition using the Dirichlet heat kernel}

\bpr \label{prop:decompose-ball-heat kernel}

For a function $f(x,t) \in C^{2,1}_{x,t}\left(B_1\times(0,\infty)\right) \cap C_{x,t}\left(\bar{B_1} \times (0,\infty)\right) \cap C_{x,t}\left(B_1 \times [0,\infty)\right)$, one can decompose it into

 \be\label{eq:heat-ball-decomposition}
  \begin{split}
f(x,t)=&\int_{B_1}\Gamma_1(x,y,t)f(y,0) dy + \int_0^t\int_{\pa B_1}E_1(x,y,t-s)f(y,s) d\sigma_y ds  \\
 &+  \int_0^t\int_{ B_1}\Gamma_1(x,y,t-s)\Big[(\pa_s - \Delta_y)f(y,s)\Big] dy ds ,
  \end{split}
\ee
  where
  \[
  E_1(x,y,t)= -\pa_{\nu_y}\Gamma_1(x,y,t) , \quad x \in B_1, y\in \pa B_1, t>0.
  \]

\epr

\bpf
 Let
 \[
 u(x,t)= \int_0^t \int_{\pa B_1} E_1(x,y,t-s) f(y,s) d\sigma_y ds,
 \]

then integration by parts and divergence theorem provide
 \be\label{eq:integration by parts-E_1}
  \begin{split}
u(x,t)=&\int_0^t \int_{B_1} - \divg_y  \Big( \nabla_y \Gamma_1(x,y,t-s)f(y,s) \Big)   dy  ds   \\
 =& \int_0^t \int_{B_1} -\Delta_y \Gamma_1(x,y,t-s)f(y,s) - \nabla_y  \Gamma_1(x,y,t-s)  \cdot \nabla_y f(y,s)  dy ds \\
  =& \int_0^t \int_{B_1} -\pa_t \Gamma_1(x,y,t-s)f(y,s) +  \Gamma_1(x,y,t-s)  \Delta_y f(y,s)  dy ds \\
    =& \int_0^t \int_{B_1} \pa_s \Big( \Gamma_1(x,y,t-s)f(y,s) \Big) - \Gamma_1(x,y,t-s) \Big[(\pa_s- \Delta_y)  f(y,s) \Big] dy ds \\
    =& \lim_{\tau \to 0^+}\int_{B_1}   \Gamma_1(x,y,\tau)f(y,t) dy  - \int_{B_1}   \Gamma_1(x,y,t)f(y,0) dy \\
    &- \int_0^t \int_{B_1} \Gamma_1(x,y,t-s) \Big[(\pa_s- \Delta_y)  f(y,s) \Big] dy ds \\
    =&F(x,t)-v(x,t)-w(x,t)
  \end{split}
\ee
where we used the equation of $\Gamma_1$ in the third line. By  \eqref{eq:dirichlet-heat-kernel-ball-pde-Gamma_1} and the Duhamel's principle, we deduce that
\be
F(x,t)=f(x,t) ,\quad  (x,t) \in B_1 \times [0,\infty),
\ee
 $v(x,t)$ satisfies
\be\label{eq:decompose-heat kernel-v-initial}
  \left\{
\begin{aligned}
 \pa_t v -\Delta v & = 0 \quad                           \mbox{in }   B_1 \times (0, \infty),\\
 v &    =0    \quad  \mbox{on } \pa B_1 \times (0, \infty) ,\\
 v & = f  \quad  \mbox{on } B_1 \times \{0\},
\end{aligned}
\right.
\ee
and $w(x,t)$ satisfies
\be\label{eq:decompose-heat kernel-w-unhomo}
  \left\{
\begin{aligned}
 \pa_t w -\Delta w & = \pa_t f -\Delta f \quad                           \mbox{in }   B_1 \times (0, \infty),\\
 w  &    =0    \quad  \mbox{on } \pa B_1 \times (0, \infty) ,\\
 w  & = 0  \quad  \mbox{on } B_1 \times \{0\}.
\end{aligned}
\right.
\ee
By the continuity of $f$ and the property of $v$ and $w$, one finds
\[
\lim_{B_1 \ni x\to x_0}u(x,t) =\lim_{B_1 \ni x\to x_0}F(x,t)=f(x_0,t), \quad x_0 \in \pa B_1, t>0.
\]
Therefore, $u(x,t)$ is a solution to
\be\label{eq:decompose-heat kernel-u-boundary}
  \left\{
\begin{aligned}
 \pa_t u -\Delta u & = 0 \quad                           \mbox{in }   B_1 \times (0, \infty),\\
 u  &    =f    \quad  \mbox{on } \pa B_1 \times (0, \infty) ,\\
 u  & = 0  \quad  \mbox{on } B_1 \times \{0\}.
\end{aligned}
\right.
\ee
with
\[
\lim_{t \to 0^+}u(x,t) =f(x,0), \quad \forall x \in \pa B_1.
\]
By the theory of parabolic equation, the decomposition \eqref{eq:heat-ball-decomposition} holds.
\epf

\brem \label{rem:interpretation-E_1-heat}
$E_1(x,y,t)$ is characterized as
\be\label{eq:dirichlet-heat-kernel-ball-pde-boundary-E_1}
  \left\{
\begin{aligned}
 \pa_t E_1 -\Delta_x E_1 & = 0 \quad                           \mbox{in }   B_1 \times (0, \infty),\\
 E_1 &    =\delta_x(y)\delta_0(t)    \quad  \mbox{on } \pa B_1 \times (0, \infty) ,\\
 E_1(0) & = 0   \quad  \mbox{on } B_1\times\{0\}  .
\end{aligned}
\right.
\ee

\erem

Under this scheme, we can solve $\phi_1^y(x,t)$ as
\[
\phi_1^y(x,t) = \int_0^t\int_{\pa B_1} E_1(x,z,t-s)\Gamma(z,y,s) d\sigma_z ds = \int_0^t\int_{\pa B_1} \Gamma(z,y,t-s)\Big(-\pa_{\nu_y}\Gamma_1(x,z,s) \Big)d\sigma_z ds  .
\]
And naturally, we derive a probabilistic representation:
\[
\int_0^t\int_{\pa B_1} E_1(x,z,t-s)\Gamma(z,y,s) d\sigma_z ds = \mathbf{E}^x \left( t>\tau_1, \Gamma(W_{\tau_1} , y,t - \tau_1)\right).
\]

\subsubsection{Related kernels}
\bpr \label{prop:soluton-u-ball-heat kernel-F_1-Gamma_1}
Let
\[
F_1(x,y,t)=\int_0^t E_1(x,y,t-s)ds, \quad x \in B_1, y \in \pa B_1, t>0,
\]
then it is characterized as
\be\label{eq:dirichlet-heat-kernel-ball-pde-boundary-F_1}
  \left\{
\begin{aligned}
 \pa_t F_1 -\Delta_x F_1 & = 0 \quad                           \mbox{in }   B_1 \times (0, \infty),\\
 \pa_t F_1 &    =0     \quad  \mbox{on } \pa B_1 \times (0, \infty) ,\\
  F_1(0)&=\delta_x(y)  \quad \mbox{on }  \pa B_1 \times \{0\},\\
 F_1(0) & = 0   \quad  \mbox{on } B_1 \times \{0\}.
\end{aligned}
\right.
\ee
and
\be \label{eq:-heat-ball-integral-F_1 and Gamma_1}
\int_{\pa B_1}  F_1(x,y,t) d \sigma_y + \int_{ B_1} \Gamma_1(x,y,t) dy = 1, \quad \forall
x \in B_1, t>0.
\ee

\epr

\bpf

Let
\[
u(x,t)= \int_{\pa B_1}  F_1(x,y,t) \phi_b(y) d \sigma_y + \int_{ B_1} \Gamma_1(x,y,t) \phi_i(y)dy ,
\]
with
\[
 \phi_b(x) \in C(\pa B_1) , \quad \phi_i(x) \in C(B_1) \cap L^1(B_1).
\]

We denote $\Phi_b(x) \in C(\bar{B_1}) \cap C^2(B_1)$ as a harmonic extension of $\phi_b(x)$, that
\be\label{eq:Laplace-ball-harmonic-extension}
  \left\{
\begin{aligned}
 -\Delta \Phi_b & =  0 \quad                           \mbox{in }   B_1 \\
 \Phi_b  &    =  \phi_b    \quad  \mbox{on } \pa B_1 .
\end{aligned}
\right.
\ee

If $x \in B_1, t>0$,
 \be
  \begin{split}
\pa_t u - \Delta u =& \int_{\pa B_1} \Big[(\pa_t - \Delta_x)F_1(x,y,t)\Big] \phi_b(y) d\sigma_y  +  \int_{ B_1} \Big[(\pa_t - \Delta_x)\Gamma_1(x,y,t)\Big] \phi_i(y)dy \\
=& \int_0^t\int_{\pa B_1} \Big[(\pa_t - \Delta_x)E_1(x,y,t-s)\Big] \phi_b(y) d\sigma_y +
 \lim_{\tau \to 0^+}\int_{\pa B_1} E_1(x,y,\tau) \phi_b(y) d\sigma_y  \\
=& 0.
 \end{split}
\ee
For $x \in \pa B_1, t>0$,
 \be    \label{eq:boundary-E_1-0 for space integral}
  \begin{split}
\pa_t u & =\int_{\pa B_1} \pa_t F_1(x,y,t) \phi_b(y) d\sigma_y  +  \int_{ B_1} \pa_t \Gamma_1(x,y,t) \phi_i(y)dy \\
&= \int_{\pa B_1} E_1(x,y,t) \phi_b(y) d\sigma_y  \\
 &=-\int_{\pa B_1} \pa_{\nu_y} \Gamma_1(x,y,t) \phi_b(y) d\sigma_y\\
 &=-\int_{ B_1} \divg_y \Big( \nabla_y \Gamma_1(x,y,t) \Phi_b(y)  \Big) dy\\
 &=-\int_{ B_1} \Delta_y \Gamma_1(x,y,t) \Phi_b(y)  dy-\int_{ B_1} \nabla_y \Gamma_1(x,y,t)   \cdot \nabla_y \Phi_b(y)  dy \\
 &=-\int_{ B_1} \pa_t \Gamma_1(x,y,t) \Phi_b(y)  dy-\int_{ B_1} \divg_y \Big( \Gamma_1(x,y,t)\nabla_y \Phi_b(y) \Big)  dy  + \int_{ B_1}  \Gamma_1(x,y,t) \Delta_y\Phi_b(y)  dy\\
  &=-\int_{ B_1} \pa_t \Gamma_1(x,y,t) \Phi_b(y)  dy-\int_{ \pa B_1}  \Gamma_1(x,y,t) \pa_{\nu_y} \Phi_b(y)   d\sigma_y  + \int_{ B_1}  \Gamma_1(x,y,t) \Delta_y\Phi_b(y)  dy\\
  &=0.
 \end{split}
\ee
Here we used the fact that
\[
\Gamma_1(x,y,t) = \pa_t \Gamma_1(x,y,t)=0, \quad (x,t) \in \pa B_1 \times(0,\infty).
\]
Noticing
\[
\lim_{t \to 0^+} \int_{\pa B_1}  F_1(x,y,t) \phi_b(y) d \sigma_y  = \lim_{t \to 0^+} \int_0^t \int_{\pa B_1}  E_1(x,y,t-s) \phi_b(y) d \sigma_y ds = \phi_b(x)
\]
for $x \in \pa B_1$. Then, together with the property of $\Gamma_1$ and $E_1$,  we find that $u(x,t)$ satisfies the initial value problem
\be\label{eq:dirichlet-heat-kernel-ball-pde-boundary-u-F_1 and Gamma_1}
  \left\{
\begin{aligned}
 \pa_t u -\Delta u & = 0 \quad                           \mbox{in }   B_1 \times (0, \infty),\\
 \pa_t u &    =0    \quad  \mbox{on } \pa B_1 \times (0, \infty) ,\\
  u(0)&=\phi_b(x)  \quad \mbox{on }  \pa B_1,\\
 u(0) & = \phi_i(x)   \quad  \mbox{on } B_1.
\end{aligned}
\right.
\ee
The proposition is proved if we take $u \equiv 1$ on $\bar{B_1} \times [0,\infty)$ as a particular solution of \eqref{eq:dirichlet-heat-kernel-ball-pde-boundary-u-F_1 and Gamma_1} with the initial datum $\phi_b(x)\equiv 1 $ and $\phi_i(x)\equiv 1 $ .

\epf

\bpr \label{prop:soluton-u-ball-heat kernel-H_1}

Let
\[
H_1(x,y,t)=\int_0^t E_1(xe^{-s},y,t-s)ds, \quad x \in B_1, y \in \pa B_1, t>0,
\]
then we have
\be\label{eq:dirichlet-heat-kernel-ball-pde-boundary-H_1}
  \left\{
\begin{aligned}
 \pa_t H_1 -\Delta_x H_1 & = \mathbf{F}(x,t) \quad                           \mbox{in }   B_1 \times (0, \infty),\\
 \pa_t H_1 + \pa_{\nu_x}   H_1&    =0     \quad  \mbox{on } \pa B_1 \times (0, \infty) ,\\
  H_1(0)&=\delta_x(y)  \quad \mbox{on }  \pa B_1,\\
   H_1(0)&=0  \quad \mbox{on }   B_1,
\end{aligned}
\right.
\ee
with
\[
|\mathbf{F}(x,t)| \to 0, \quad  t \to 0^+.
\]

\epr

\bpf
Let
\[
u(x,t) = \int_{\pa B_1 } H_1(x,y,t) \phi_b(y) d\sigma_y
\]
and $\Phi_b(x) \in C(\bar{B_1}) \cap C^2(B_1)$ be the harmonic extension of $\phi_b(x)$, satisfying \eqref{eq:Laplace-ball-harmonic-extension}. Then, integration by parts gives
\be\label{eq:integration by parts-H_1}
  \begin{split}
u(x,t)=& -\int_0^t \int_{ B_1}  \Delta_y \Gamma_1(xe^{-s},y,t-s)\Phi_b(y)  dy  ds        \\
&- \int_0^t \int_{ B_1} \Gamma_1(xe^{-s},y,t-s) (-\Delta_y \Phi_b(y)) dyds \\
=&-\int_0^t \int_{ B_1} \pa_t   \Big(     \Gamma_1(xe^{-s},y,t-s)\Phi_b(y)       \Big) dy ds \\
=&\int_0^t \int_{ B_1} \pa_s   \Big(     \Gamma_1(xe^{-s},y,t-s)\Phi_b(y)       \Big) dy ds  +  \int_0^t \int_{ B_1} e^{-s} \pa_{e_x} \Gamma_1(xe^{-s},y,t-s)\Phi_b(y)      dyds   \\
=&  \lim_{\tau \to 0^+}\int_{ B_1}\Gamma_1(xe^{-t},y,\tau)\Phi_b(y)       dy - \int_{ B_1}\Gamma_1(x,y,t)\Phi_b(y) dy \\
&+\int_0^t \int_{ B_1} e^{-s} \pa_{e_x} \Gamma_1(xe^{-s},y,t-s)\Phi_b(y)      dyds\\
=&   \Phi_b(xe^{-t})       - v(x,t) +\int_0^t \int_{ B_1} e^{-s} \pa_{e_x} \Gamma_1(xe^{-s},y,t-s)\Phi_b(y)      dyds  ,
  \end{split}
\ee
where $v(x,t)=\int_{ B_1}\Gamma_1(x,y,t)\Phi_b(y) dy$ satisfies
\be\label{eq:decompose-heat-kernel-H_1-v-initial}
  \left\{
\begin{aligned}
 \pa_t v -\Delta v & = 0 \quad                           \mbox{in }   B_1 \times (0, \infty),\\
 v &    =0    \quad  \mbox{on } \pa B_1 \times (0, \infty) ,\\
 v & = \Phi_b  \quad  \mbox{on } B_1 \times \{0\},
\end{aligned}
\right.
\ee
and $e_x$ is the unit vector along the $x$ direction, i.e.
\[
e_x = \frac{x}{|x|}, \quad x \neq 0.
\]
We use the symmetry of $\Gamma_1(x,y,t)$ with respect to $x,y$, and the property of $E_1$ in  Remark \ref{rem:interpretation-E_1-heat}, to know that
\be \label{eq:initial-pa_ex-Gamma_1-boundary}
\lim_{t\to 0^+} \int_{B_1} -\pa_{e_x} \Gamma_1(x,y,t)\Phi_b(y) dy=0 , \quad x \in \pa B_1,
\ee
getting
\[
\lim_{t\to 0^+} u(x,t) = \phi_b(x) , \quad x \in \pa B_1.
\]
By the property of $\Gamma_1$, that $\Gamma_1$ behaves like $\Gamma$ for $t$ is  close to 0 and $x,y$ are away from the boundary $\pa B_1$, we find
\be \label{eq:initial-pa_ex-Gamma_1-interior}
\lim_{t\to 0^+}  \int_0^t \int_{B_{1-\epsilon}} -\pa_{e_x} \Gamma_1(xe^{-s},y,t-s)\Phi_b(y) dy ds=0 , \quad x \in \bar{B_1} - \{0\}, 0<\epsilon<1
\ee
getting
\[
\lim_{t\to 0^+} u(x,t) = 0 ,
\]
for $x \in B_1$ almost everywhere.

Next, if $x \in B_1, t>0$,
\be
  \begin{split}
&\pa_t u - \Delta u \\
=& \int_{\pa B_1} \Big[(\pa_t - \Delta_x)H_1(x,y,t)\Big] \phi_b(y) d\sigma_y   \\
=& \int_0^t\int_{\pa B_1} \Big[(1-e^{-2s})\pa_t E_1(xe^{-s},y,t-s)\Big] \phi_b(y) d\sigma_y  ds  +
 \lim_{\tau \to 0^+}\int_{\pa B_1} E_1(xe^{-t},y,\tau) \phi_b(y) d\sigma_y  \\
 =&\int_0^t   (1-e^{-2s})ds \int_{\pa B_1}    \pa_t E_1(xe^{-s},y,t-s)  \phi_b(y) d\sigma_y \\
 =&\mathbf{F}(x,t) ,   \quad |\mathbf{F}(x,t)| \to 0, \quad \mbox{as } \quad  t \to 0^+.
 \end{split}
\ee
For $x \in \pa B_1, t>0$,
 \be
  \begin{split}
\pa_t u + \pa_{\nu_x} u & =\int_{\pa B_1} \Big[(\pa_t + \pa_{e_x} ) H_1(x,y,t)\Big] \phi_b(y) d\sigma_y   \\
 &=  \int_0^t  \int_{\pa B_1}  \Big( \pa_t E_1(xe^{-s},y,t-s)  + e^{-s}\pa_{e_x}E_1(xe^{-s},y,t-s) \Big)   \phi_b(y) d\sigma_y ds\\
 & \quad + \lim_{\tau \to 0^+}\int_{\pa B_1} E_1(xe^{-t},y,\tau) \phi_b(y) d\sigma_y \\
 &=  \int_0^t  \int_{\pa B_1} -\pa_s \Big( E_1(xe^{-s},y,t-s) \phi_b(y)  \Big)    d\sigma_y ds + \lim_{\tau \to 0^+}\int_{\pa B_1} E_1(xe^{-t},y,\tau) \phi_b(y) d\sigma_y \\
   &=  -\lim_{\tau \to 0^+}\int_{\pa B_1} E_1(xe^{-t},y,\tau) \phi_b(y) d\sigma_y  + \lim_{\tau \to 0^+}\int_{\pa B_1} E_1(xe^{-t},y,\tau) \phi_b(y) d\sigma_y   \\
   & \quad +\int_{\pa B_1} E_1(x,y,t) \phi_b(y) d\sigma_y \\
   &= \int_{\pa B_1} E_1(x,y,t) \phi_b(y) d\sigma_y =0,
 \end{split}
\ee
where we used \eqref{eq:boundary-E_1-0 for space integral} in the last line.

\epf

\newpage

\section{Green's function for heat equation with a dynamical boundary condition in the unit ball}\label{sec:green}

In this section, we are aim to find the Green's function $G_1(x,y,t)$ to \eqref{eq:heat--dynamical} for $\Omega = B_1$, such that

\be\label{eq:heat--dynamical-G_1}
  \left\{
\begin{aligned}
\pa_t G_1 -\Delta_x G_1 & = 0 \quad                  \mbox{in }   B_1 \times (0,\infty),\\
 \pa_t G_1+\pa_{\nu_x} G_1  &    =0    \quad  \mbox{on } \pa B_1\times (0,\infty) ,
\end{aligned}
\right.
\ee
and
\[
   u(x,t) = \int_{B_1} G_1(x,y,t) \phi_i(y) dy + \int_{\partial B_1} G_1(x,y,t) \phi_b(y) d\sigma_y
   \]
is a solution to the initial value problem \eqref{eq:heat--dynamical} for $\Omega = B_1$.
And \(G_1(x,y,t)\) should satisfy
\be \label{eq:-heat-ball-integral-G_1}
\int_{\pa B_1}  G_1(x,y,t) d \sigma_y + \int_{ B_1} G_1(x,y,t) dy = 1, \quad \forall
x \in B_1, t>0,
\ee
\[
 G_1(x,y,t) = G_1(y,x,t) , \quad x,y \in B_1, t>0.
\]
Here, we introduce an implicit formula of $G_1$ in an analytic view. Follow the idea that we present $\Gamma_1(x,y,t)$ by the eigenfunctions and eigenvalues to problem  \eqref{eq:dirichlet-laplace-ball-eigen}. We can present such fundamental solution $G_1(x,y,t)$ to the given initial-boundary value problem with a dynamic boundary condition, exploiting an eigenfunction expansion tailored to the specific boundary conditions.

Let  $(\lambda_i , \phi_i )$ be the eigen pair for the eigenvalue problem
\be\label{eq:Neumann-laplace-ball-eigen}
  \left\{
\begin{aligned}
 -\Delta \psi_i & = \lambda_i \psi_i \quad                           \mbox{in }   B_1 ,\\
 \pa_\nu \psi_i & = \lambda_i \psi_i \quad  \mbox{on } \pa B_1,
\end{aligned}
\right.
\ee
then the kernel  \(G_1(x,y,t)\) is given by
   \[
   G_1(x,y,t) = \sum_{k} e^{-\lambda_k t} \psi_k(x) \psi_k(y).
   \]
This allows the solution of \eqref{eq:heat--dynamical} to be written as a series:
   \[
   u(x,t) = \sum_{k} e^{-\lambda_k t} \left( \int_{B_1} \phi_i(y) \psi_k(y) dy + \int_{\partial B_1} \phi_b(y) \psi_k(y) d\sigma_y \right) \psi_k(x),
   \]
   where the eigenfunctions \(\psi_k\) are orthogonal with respect to the inner product over \(B_1 \cup \partial B_1\).

The elliptic problem \eqref{eq:Neumann-laplace-ball-eigen} is of importance in research. And yet, there is no detailed study on this problem.

In view of recent research results, I believe there is also an probability representation of $G_1$ that would provide more accurate knowledge about this kernel and equation \eqref{eq:heat--dynamical}. We shall look into this probabilistic direction in the future research.

In the  next section, we use the similar method that was used in the half space, building some approximations to $G_1$ and the solution $u$ of \eqref{eq:heat--dynamical}.

\newpage

\section{Approximations for heat equation with a dynamical boundary condition in the unit ball}\label{sec:appro}

In this section, we apply and improve the methods used for the case that $\om = \R^n_+$ and $\om= \R^n - B_1$ to build approximations for Green's function and solution of the unit ball case for \eqref{eq:heat--dynamical}.

\subsection{Approximation of the Green's function $G_1$} \label{subsec:appro to green}
We use $\Gamma_1$ to build an approximation to $G_1$. Sadly, it possesses a singularity at the origin because $\pa_\nu$ is not of constant coefficients on $\pa B_1$. For the case of upper half space, such problem is missing because $\pa_\nu=-\pa_{x_n}$ on $\pa \R^n_+$.
\bthm \label{thm:G-1 approximation}

Let

\[
\mathcal{G}_1(x,y,t)=\Gamma_1(x,y,t)+\mathcal{H}_1(x,y,t), \quad  x \in B_1 - \{0\}, y \in B_1, t>0
\]
where
\[
\mathcal{H}_1(x,y,t)= -\int_0^t \pa_{e_x}\Gamma_1(xe^{-s},y,t-s)ds, \quad  x \in B_1 - \{0\}, y \in B_1, t>0.
\]

For
\[
u(x,t)= \int_{\pa B_1}  \mathcal{G}_1(x,y,t) \phi_b(y) d \sigma_y + \int_{ B_1} \mathcal{G}_1(x,y,t) \phi_i(y)dy ,
\]
with
\[
 \phi_b(x) \in C(\pa B_1) , \quad \phi_i(x) \in C(B_1) \cap L^1(B_1),
\]
and $\Phi_b(x) \in C(\bar{B_1}) \cap C^2(B_1)$ be the harmonic extension of $\phi_b(x)$ as we previous done.
Then $u$ is a solution to the initial value problem
\be\label{eq:heat--ball-pde-dynamical-u-O(t)}
  \left\{
\begin{aligned}
 \pa_t u -\Delta_x u & = \mathbf{F}(x,t) \quad                           \mbox{in }  \left( B_1-\{0\} \right) \times (0, \infty),\\
 \pa_t u+ \pa_{\nu_x}   u&    =0     \quad  \mbox{on } \pa B_1 \times (0, \infty) ,\\
  u(0)&=\phi_b  \quad \mbox{on }  \pa B_1,\\
   u(0)&=\phi_i  \quad \mbox{on }   B_1,
\end{aligned}
\right.
\ee
with
\[
|\mathbf{F}(x,t)| \to 0, \quad  t \to 0^+.
\]

\ethm

\newpage

\bpf

By \eqref{eq:initial-pa_ex-Gamma_1-interior}, we find
\[
\lim_{t \to 0^+} \int_{ B_1} \mathcal{H}_1(x,y,t) \phi_i(y)dy =0, \quad x \in  \bar{B_1}-\{0\}.
\]
By the fact that $\Gamma_1(x,y,t) =0 , y \in \pa B_1$,
we find
\[
\lim_{\tau \to 0^+} \int_{\pa B_1} -\pa_{e_x}\Gamma_1(x^{-s},y,\tau) \phi_b(y) d\sigma_y = 0,\quad x \in  B_1-\{0\},s>0,
\]
therefore,
\[
\lim_{t \to 0^+} \int_{  \pa B_1} \mathcal{H}_1(x,y,t) \phi_b(y)d\sigma_y = 0 ,\quad x \in  B_1-\{0\} .
\]
From the property of $E_1(x,y,t)$ and the symmetry of $\Gamma_1(x,y,t)$, we characterize that
\[
\lim_{t \to 0^+} -\pa_{\nu_x}\Gamma_1(x,y,t) = \delta_x(y) ,\quad x,y \in \pa B_1,
\]
therefore,
\[
\lim_{t \to 0^+} \int_{  \pa B_1} \mathcal{H}_1(x,y,t) \phi_b(y)d\sigma_y = \phi_b(x) ,\quad x \in \pa B_1 .
\]
Conclude above equalities, we have
\[
\lim_{t \to 0^+}u(x,t)= \phi_i(x), \  x\in  B_1-\{0\}, \quad \lim_{t \to 0^+}u(x,t)= \phi_b(x), \  x\in   \pa B_1.
\]

For $x \in \pa B_1, t>0$,
\be
\begin{split}
&(\pa_t + \pa_{\nu_x})\mathcal{H}_1(x,y,t)\\
&= \lim_{\tau \to 0^+}-\pa_{e_x}\Gamma_1(xe^{-t},y,\tau) \\
&+ \int_0^t -\pa_t\pa_{e_x} \Gamma_1(xe^{-s},y,t-s)ds+\int_0^t -e^{-s}\pa_{e_x}\pa_{e_x} \Gamma_1(xe^{-s},y,t-s)ds\\
&= \lim_{\tau \to 0^+}-\pa_{e_x}\Gamma_1(xe^{-t},y,\tau) + \int_0^t \pa_s  \Big(\pa_{e_x} \Gamma_1(xe^{-s},y,t-s)   \Big)   ds\\
&= \pa_{e_x} \Gamma_1(x,y,t).
\end{split}
\ee

Consequently,
\[
(\pa_t + \pa_{\nu_x})\mathcal{G}_1(x,y,t) = 0, \mbox{ on } \pa B_1 \times (0,\infty).
\]

 \be
  \begin{split}
&\pa_t u - \Delta u = \int_{ B_1} \Big[(\pa_t - \Delta_x)\mathcal{H}_1(x,y,t)\Big] \phi_i(y) dy + \int_{\pa B_1} \Big[(\pa_t - \Delta_x)\mathcal{H}_1(x,y,t)\Big] \phi_b(y) d\sigma_y   \\
&= \lim_{\tau \to 0^+}\int_{ B_1} -\pa_{e_x}\Gamma_1(xe^{-t},y,\tau) \phi_i(y) dy  +
 \lim_{\tau \to 0^+}\int_{\pa B_1} -\pa_{e_x}\Gamma_1(xe^{-t},y,\tau) \phi_b(y) d\sigma_y  \\
 &  +  \int_0^t\int_{ B_1} \Big[(e^{-2s}-1)\pa_t \pa_{e_x}\Gamma_1(xe^{-s},y,t-s)\Big] \phi_i(y) dy  ds \\
 &  +  \int_0^t\int_{\pa B_1} \Big[(e^{-2s}-1)\pa_t \pa_{e_x}\Gamma_1(xe^{-s},y,t-s)\Big] \phi_b(y) d\sigma_y  ds \\
 &  - \int_0^t\int_{ B_1} \frac{2}{|x|} e^{-s}\Delta_x   \Gamma_1(xe^{-s},y,t-s)   \phi_i(y) dy  ds  - \int_0^t\int_{ \pa B_1} \frac{2}{|x|} e^{-s}\Delta_x   \Gamma_1(xe^{-s},y,t-s)   \phi_b(y) d\sigma_y  ds\\
 &  + \int_0^t\int_{ B_1} \frac{n-1}{|x|^2} \pa_{e_x}  \Gamma_1(xe^{-s},y,t-s)   \phi_i(y) dy  ds   + \int_0^t\int_{ \pa B_1} \frac{n-1}{|x|^2} \pa_{e_x}   \Gamma_1(xe^{-s},y,t-s)   \phi_b(y) d\sigma_y  ds\\
  &  + \int_0^t\int_{ B_1} \frac{2}{|x|}   e^{-s}   \pa_{e_x} \pa_{e_x} \Gamma_1(xe^{-s},y,t-s)   \phi_i(y) dy  ds   + \int_0^t\int_{ \pa B_1} \frac{2}{|x|}   e^{-s}   \pa_{e_x} \pa_{e_x}  \Gamma_1(xe^{-s},y,t-s)   \phi_b(y) d\sigma_y  ds\\
& =    \int_0^t(e^{-2s}-1)ds \left( \int_{ B_1} \pa_t \pa_{e_x}\Gamma_1(xe^{-s},y,t-s)  \phi_i(y) dy    +  \int_{\pa B_1} \pa_t \pa_{e_x}\Gamma_1(xe^{-s},y,t-s) \phi_b(y) d\sigma_y  \right)  \\
 &  - \frac{2}{|x|} \int_0^t  e^{-s}  ds \left( \int_{ B_1}  \pa_t \Gamma_1(xe^{-s},y,t-s)   \phi_i(y) dy  +  \int_{ \pa B_1}  \pa_t \Gamma_1(xe^{-s},y,t-s)   \phi_b(y) d\sigma_y  \right) \\
 &  + \frac{n-1}{|x|^2}\int_0^t ds \left( \int_{ B_1}  \pa_{e_x}  \Gamma_1(xe^{-s},y,t-s)   \phi_i(y) dy    + \int_{ \pa B_1}  \pa_{e_x}   \Gamma_1(xe^{-s},y,t-s)   \phi_b(y) d\sigma_y     \right)\\
  &  +  \frac{2}{|x|}   \int_0^te^{-s} ds  \left( \int_{ B_1}   \pa_{e_x} \pa_{e_x} \Gamma_1(xe^{-s},y,t-s)   \phi_i(y) dy  + \int_{ \pa B_1}  \pa_{e_x} \pa_{e_x}  \Gamma_1(xe^{-s},y,t-s)   \phi_b(y) d\sigma_y  \right)\\
  & =    \int_0^t(e^{-2s}-1)ds \left( \int_{ B_1} \pa_t \pa_{e_x}\Gamma_1(xe^{-s},y,t-s)  \phi_i(y) dy    +  \int_{\pa B_1} \pa_t \pa_{e_x}\Gamma_1(xe^{-s},y,t-s) \phi_b(y) d\sigma_y  \right)  \\
 &  - \frac{2}{|x|} \int_0^t  e^{-s}  ds \left( \int_{ B_1}  \pa_t \Gamma_1(xe^{-s},y,t-s)   \phi_i(y) dy  +  \int_{ \pa B_1}  \pa_t \Gamma_1(xe^{-s},y,t-s)   \phi_b(y) d\sigma_y  \right) \\
 &  - \frac{n-1}{|x|^2} \left( \int_{ B_1}  \mathcal{H}_1(x,y,t)   \phi_i(y) dy    + \int_{ \pa B_1}  \mathcal{H}_1(x,y,t)   \phi_b(y) d\sigma_y     \right)\\
  &  +  \frac{2}{|x|}   \int_0^te^{-s} ds  \left( \int_{ B_1}   \pa_{e_x} \pa_{e_x} \Gamma_1(xe^{-s},y,t-s)   \phi_i(y) dy  + \int_{ \pa B_1}  \pa_{e_x} \pa_{e_x}  \Gamma_1(xe^{-s},y,t-s)   \phi_b(y) d\sigma_y  \right)\\
   &  =    \int_0^t(e^{-2s}-1)ds \left( \int_{ B_1} \pa_t \pa_{e_x}\Gamma_1(xe^{-s},y,t-s)  \phi_i(y) dy    +  \int_{\pa B_1} \pa_t \pa_{e_x}\Gamma_1(xe^{-s},y,t-s) \phi_b(y) d\sigma_y  \right)  \\
 &  - \frac{2}{|x|} \int_0^t  e^{-s}  ds \left( \int_{ B_1}  \pa_t \Gamma_1(xe^{-s},y,t-s)   \phi_i(y) dy  +  \int_{ \pa B_1}  \pa_t \Gamma_1(xe^{-s},y,t-s)   \phi_b(y) d\sigma_y  \right) \\
 &  - \frac{n-1}{|x|^2} \Big( u(x,t)-v(x,t)    \Big)\\
  &  +  \frac{2}{|x|}   \int_0^te^{-s} ds  \left( \int_{ B_1}   \pa_{e_x} \pa_{e_x} \Gamma_1(xe^{-s},y,t-s)   \phi_i(y) dy  + \int_{ \pa B_1}  \pa_{e_x} \pa_{e_x}  \Gamma_1(xe^{-s},y,t-s)   \phi_b(y) d\sigma_y  \right)\\
  &  =\mathbf{F}(x,t),  \quad  x \in B_1-{0}, t>0 .
 \end{split}
\ee
where
\[
v(x,t) = \int_{\pa B_1}  \Gamma_1(x,y,t) \phi_b(y) d \sigma_y + \int_{ B_1} \Gamma_1(x,y,t) \phi_i(y)dy = \int_{ B_1} \Gamma_1(x,y,t) \phi_i(y)dy
\]
is a solution to the heat equation with Dirichlet boundary condition given by
\be\label{eq:heat kernel-v-initial-phi_i}
  \left\{
\begin{aligned}
 \pa_t v -\Delta v & = 0 \quad                           \mbox{in }   B_1 \times (0, \infty),\\
 v &    =0    \quad  \mbox{on } \pa B_1 \times (0, \infty) ,\\
 v & = \phi_i \quad  \mbox{on } B_1 \times \{0\}.
\end{aligned}
\right.
\ee
So,
\[
\lim_{t \to 0^+}(u(x,t)-v(x,t))=0, \quad     x \in B_1-{0}, t>0 .
\]
Consequently,
\[
 |\mathbf{F}(x,t)| \to 0, \quad \mbox{as } \quad  t \to 0^+ ,
 \]
and $u(x,t)$ is a solution to \eqref{eq:heat--ball-pde-dynamical-u-O(t)}.

\epf

\subsection{Approximation of the solution $u$ }\label{subsec:appro to solution}

Here, we build an approximation to the solution of problem \eqref{eq:heat--dynamical} for $\om = B_1$.

\bthm \label{thm:u approximation}

Let

\[
\tilde{\Gamma_1 }(x,y,t)=\Gamma_1(x,y,t)-\int_0^t\int_{\pa B_1}K_1(x,z,t-s)\pa_{\nu_z}\Gamma_1(z,y,s) d\sigma_z ds , \quad  x \in B_1 , y \in B_1, t>0
\]
and
\be
\begin{split}
\tilde{H}_1(x,y,t)&= K_1(x,y,t) - \int_{B_1} \Gamma_1(x,z,t)P_1(z,y)dz \\
+&\int_0^t \int_{\pa B_1} \int_{B_1} K_1(x,w,t-s)\pa_{\nu_w}\Gamma_1(w,z,s)P_1(z,y)dzd\sigma_wds, \quad  x \in B_1 , y  \in \pa  B_1, t>0,
\end{split}
\ee
where $K_1(x,y,t)=P_1(xe^{-t},y)$ is the Green's function to \eqref{eq:Laplace--dynamical} for $\om = B_1$, and $P_1(x,y)$ is the Poisson kernel of the unit ball presented in formula \eqref{eq:poisson kernel-ball}.

For
\[
u(x,t)= \int_{\pa B_1} \tilde{H}_1(x,y,t) \phi_b(y) d \sigma_y + \int_{ B_1} \tilde{\Gamma_1 }(x,y,t) \phi_i(y)dy ,
\]
with
\[
 \phi_b(x) \in C(\pa B_1) , \quad \phi_i(x) \in C(B_1) \cap L^1(B_1),
\]
and $\Phi_b(x) \in C(\bar{B_1}) \cap C^2(B_1)$ be the harmonic extension of $\phi_b(x)$. Then $u$ is a solution to the initial value problem
\be\label{eq:heat--ball-pde-dynamical-u-O(t)}
  \left\{
\begin{aligned}
 \pa_t u -\Delta u & = \mathbf{\tilde{F}}(x,t)+\mathbf{\tilde{G}}(x,t) \quad                           \mbox{in }   B_1 \times (0, \infty),\\
 \pa_t u+ \pa_\nu   u&    =0     \quad  \mbox{on } \pa B_1 \times (0, \infty) ,\\
  u(0)&=\phi_b  \quad \mbox{on }  \pa B_1,\\
   u(0)&=\phi_i  \quad \mbox{on }   B_1,
\end{aligned}
\right.
\ee
with
\[
|\mathbf{\tilde{F}}(x,t)| \to 0, \   |x| \to 0^+, \quad  |\mathbf{\tilde{G}}(x,t)| \to 0, \   t \to 0^+.
\]

\ethm

\bpf
Let
\[
v(x,t)=\int_{B_1} \Gamma_1(x,y,t)\Phi_i(y)dy
\]
where
\[
\Phi_i(x)= \phi(x)-\Phi_b(x)=\phi(x)-\int_{\pa B_1}P_1(x,y)\phi_b(y)d\sigma_y, \quad x \in B_1.
\]
Then $v$ satisfies
\be\label{eq:heat kernel-v-initial-Phi_i}
  \left\{
\begin{aligned}
 \pa_t v -\Delta v & = 0 \quad                           \mbox{in }   B_1 \times (0, \infty),\\
 v &    =0    \quad  \mbox{on } \pa B_1 \times (0, \infty) ,\\
 v & = \Phi_i \quad  \mbox{on } B_1 \times \{0\}.
\end{aligned}
\right.
\ee
Let
\be
\begin{split}
w(x,t)&= \int_{\pa B_1}K_1(x,y,t)\phi_b(y)d\sigma_y  + \int_0^t \int_{\pa B_1}  K_1(x,y,t-s)\big( - \pa_{\nu_y}v(y,s) \big) d\sigma_y ds,
\end{split}
\ee
then $u(x,t)=v(x,y)+w(x,t)$ and $w$ satisfies
\be\label{eq:heat kernel-w-initial-Phi_b- pa_nu v}
  \left\{
\begin{aligned}
 -\Delta w & = 0 \quad                           \mbox{in }   B_1 \times (0, \infty),\\
 \pa_t w + \pa_{\nu}w &    =-\pa_\nu v    \quad  \mbox{on } \pa B_1 \times (0, \infty) ,\\
 w & = \Phi_b \quad  \mbox{on } B_1 \times \{0\},\\
 w & = \phi_b \quad  \mbox{on } \pa B_1 \times \{0\}.
\end{aligned}
\right.
\ee
Therefore $u$
satisfies
\be\label{eq:heat kernel-u-initial-phi_b-phi_i- pa_t w}
  \left\{
\begin{aligned}
 \pa_t u -\Delta u & = \pa_t w \quad                           \mbox{in }   B_1 \times (0, \infty),\\
 \pa_t u + \pa_{\nu}u &    = 0    \quad  \mbox{on } \pa B_1 \times (0, \infty) ,\\
 u & = \phi_i \quad  \mbox{on } B_1 \times \{0\},\\
 u & = \phi_b \quad  \mbox{on } \pa B_1 \times \{0\}.
\end{aligned}
\right.
\ee
From
\[
K_1(x,y,t)= c_n \frac{1-|xe^{-t}|^2}{|xe^{-t}-y|^n}, \quad x\in B_1, y\in \pa B_1, t>0,
\]
one computes
\[
\partial_t K(x,y,t) = \frac{2|x|^2 e^{-2t} |x e^{-t} - y|^2 + n(1 - |x|^2 e^{-2t})(|x|^2 e^{-2t} - x \cdot y e^{-t})}{|x e^{-t} - y|^{n+2}}.
\]
Therefore
\[
|\pa_t K_1| \to 0, \quad |x| \to 0^+.
\]
\be
\begin{split}
\pa_t w(x,t)&= \int_{\pa B_1} \pa_t K_1(x,y,t)\phi_b(y)d\sigma_y  + \int_0^t \int_{\pa B_1} \pa_t  K_1(x,y,t-s)\big( - \pa_{\nu_y}v(y,s) \big) d\sigma_y ds\\
+& \int_{\pa B_1}P_1(x,y)  \Big( - \pa_{\nu_y}v(y,t) \Big) d\sigma_y\\
&=\mathbf{\tilde{F}}(x,t) + \mathbf{\tilde{G}}(x,t),
\end{split}
\ee
with
\be
\begin{split}
\mathbf{\tilde{G}}(x,t)&= \int_{\pa B_1}P_1(x,y)  \Big( - \pa_{\nu_y}v(y,t) \Big) d\sigma_y \\
&=     \int_{\pa B_1}P_1(x,y)  \Big( - \int_{B_1} \pa_{\nu_y} \Gamma_1(y,z,t)\Phi_i(z)dz \Big) d\sigma_y  \\
&=   \int_{B_1} \Phi_i(z)dz \int_{\pa B_1}- P_1(x,y)  \pa_{\nu_y} \Gamma_1(y,z,t) d\sigma_y .
\end{split}
\ee
and
\[
|\mathbf{\tilde{F}}(x,t)| \to 0, \   |x| \to 0^+.
\]
From the property of $E_1$, one can get
\[
\lim_{t \to 0^+ } \int_{\pa B_1}- P_1(x,y)  \pa_{\nu_y} \Gamma_1(y,z,t) d\sigma_y  =0 ,\quad x \in B_1.
\]
Therefore,
\[
|\mathbf{\tilde{G}}(x,t)| \to 0, \   t \to 0^+, \ x \in B_1.
\]

\epf

\bigskip

\noindent\textbf{Declarations of interest}: none.

\bigskip


\bigskip

\noindent X. Yang

\noindent Email: \textsf{2768620240@qq.com}

\end{document}